\documentclass[leqno]{amsart}
\usepackage{amsmath}
\usepackage{mathtools}
\usepackage{amssymb}
\usepackage{amsthm}
\usepackage{relsize}
\usepackage{braket}
\usepackage{graphicx}
\usepackage{enumerate}
\usepackage[mathscr]{eucal}
\usepackage{upgreek}
\usepackage{xcolor}
\theoremstyle{plain}
\newtheorem{theorem}{Theorem}[section]
\newtheorem{problem}[theorem]{Problem}
\newtheorem{prop}[theorem]{Proposition}
\newtheorem{cor}{Corollary}[theorem]

\theoremstyle{definition}

\newtheorem{remark}{Remark}[section]

\newtheorem{example}{Example}[theorem]

\makeatother

\usepackage{hyperref}
\hypersetup{
    linkcolor=blue,
    filecolor=red,      
    urlcolor=magenta,
}

\usepackage[pagewise]{lineno}
\newcommand\restr[2]{{
  \left.\kern-\nulldelimiterspace 
  #1 
  \vphantom{\big|} 
  \right|_{#2} 
  }}

\begin{document}

\title[On best approximations in Banach spaces]{On best approximations in Banach spaces from the perspective of orthogonality}
\author{Debmalya Sain, Saikat Roy}

\address{(Sain)~Department of Mathematics, Indian Institute of Science, Bengaluru 560012, Karnataka, India}
\email{saindebmalya@gmail.com}

\address{(Roy)~Department of Mathematics, National Institute of Durgapur, West Bengal, India}
\email{saikatroy.cu@gmail.com}


\thanks{}




\date{}

\newcommand{\acr}{\newline\indent}


\thanks{The research of Dr. Debmalya Sain is sponsored by DST-SERB N-PDF Fellowship under the mentorship of Professor Apoorva Khare. Dr. Sain feels grateful to have the opportunity to acknowledge the colossal contribution of Mr. Swapan Bandyopahyay, an extraordinarily devoted teacher, towards rightly shaping the philosophies of so many young learners like himself. The research of Mr. Saikat Roy is supported by CSIR MHRD in form of Senior Research Fellowship under the supervision of Professor Satya Bagchi.} 

\begin{abstract}
We study best approximations in Banach spaces via Birkhoff-James orthogonality of functionals. To exhibit the usefulness of Birkhoff-James orthogonality techniques in the study of best approximation problems, some algorithms and distance formulae are presented. As an application of our study, we obtain some crucial inequalities, which also strengthen the classical H\"{o}lder's inequality. The relevance of the algorithms and the inequalities are discussed through concrete examples.
\end{abstract}

\maketitle

\section{Introduction}
The purpose of the present article is to study best approximations in Banach spaces, from the perspective of Birkhoff-James orthogonality of linear functionals. Recently, such a study has been carried out in the context of smooth, strictly convex, reflexive Banach spaces in \cite{Sai}. The current work strengthens those ideas in further detail and in a more general setting. To demonstrate the applicability of the results developed in this article, we obtain some distance formulae in certain special cases, which also give rise to some important classes of inequalities.\\

The symbol $\mathbb{X}$ denotes a Banach space. Unless otherwise specified, we work only with \emph{real} Banach spaces. Let $\theta$ denote the zero vector of any vector space, other than the scalar field $ \mathbb{R}. $ Let $B_{\mathbb{X}}$ and $S_{\mathbb{X}}$ denote the closed unit ball and the unit sphere of $\mathbb{X}$, respectively. We denote the collection of all extreme points of $B_\mathbb{X}$ by $Ext(B_\mathbb{X})$. Recall that $ \mathbb{X} $ is said to be \emph{strictly convex} if $ Ext(B_\mathbb{X}) = S_{\mathbb{X}}. $ The topological dual of $\mathbb{X}$ is denoted by $\mathbb{X}^*$. Note that the Banach space $\mathbb{X}$ can always be embedded into $\mathbb{X}^{**}$ via the canonical isometric isomorphism $\psi.$ Given any $f\in \mathbb{X}^*$, the norm attainment set of $f,$ denoted by $M_f,$ is defined by
$$M_f :=\{x\in S_{\mathbb{X}} : |f(x)|=\|f\|\}.$$
Let $\mathbb{L}(\mathbb{X})$ denote the collection of all bounded linear operators on $\mathbb{X}$, endowed with the usual operator norm. For any linear operator $T\in \mathbb{L}(\mathbb{X})$, the range of $T$ and the kernel of $T$ are denoted by $\mathcal{R}(T)$ and $\mathcal{N}(T)$, respectively. In similar spirit, we denote the kernel of any $f\in \mathbb{X}^*$ by $\mathcal{N}(f).$ Given any natural number $m$, let $\mathbb{H}$ denote the Hilbert space $\mathbb{R}^m$, equipped with the usual \emph{dot product} $\langle ~, ~\rangle.$ Members of $\mathbb{L}(\mathbb{H})$ are identified as matrices in the usual way. Given any $T\in \mathbb{L}(\mathbb{H})$, let $[T]$ denote the matrix representation of $T$ with respect to the standard ordered basis of $\mathbb{H}.$ Let $T^*$ denote the Hilbert adjoint of $T.$ Evidently, $[T^*]=[T]^t,$ where $[T]^t$ denotes the transpose of the matrix $[T].$\\

Approximation theory is an extensive field of research due to its diversified applications in many branches of Science. Given any element $x\in \mathbb{X}$ and a subspace $\mathbb{Y}$ of $\mathbb{X}$, distance between $x$ and $\mathbb{Y},$ denoted by $dist(x,\mathbb{Y}),$ is defined by 
$dist(x,\mathbb{Y}):=\inf \{\|x-y\|: y\in \mathbb{Y}\}.$
An element $y_0\in \mathbb{Y}$ is said to be a best approximation to $x$ out of $\mathbb{Y}$ if $dist(x,\mathbb{Y})=\|x-y_0\|.$ The existence and the uniqueness of the best approximation cannot be guaranteed, in general. However, the existence of best approximation(s) is evidently assured for finite-dimensional subspaces. Moreover, in case of a strictly convex Banach space, the best approximation is unique, whenever it exists.\\

Birkhoff-James orthogonality is of essential importance in understanding the geometry of a Banach space \cite{BS,B,J,Ja,Jb}. Given any two elements $x,y\in \mathbb{X}$, $x$ is said to be Birkhoff-James \cite{B,J,Ja} orthogonal to $y$, written as $x\perp_B y,$ if $ \|x+\lambda y\| \geq \|x\|~ \forall~\lambda \in \mathbb{R}. $
It is not difficult to see that $y_0\in \mathbb{Y}$ is a best approximation to $x$ out of $\mathbb{Y}$ if and only if $(x-y_0)\perp_B \mathbb{Y}.$ Given any non-zero element $x\in \mathbb{X}$, a member $f\in \mathbb{X}^*$ is called a support functional of $B_\mathbb{X}(\theta,\|x\|)$ at $x$, if $\|f\|=1$ and $f(x)=\|x\|.$ The point $x$ is called \emph{smooth} if the support functional of $B_\mathbb{X}(\theta,\|x\|)$ at $x$ is unique. Bhatia and \v{S}emrl completely characterized Birkhoff-James orthogonality of matrices in \cite{BS}. Based on this rudimentary result, the authors provided some distance formulae in the same article. One may consult \cite{Ga,LR} for a study of best approximations and orthogonality of matrices. We refer the readers to \cite{PSJ,S,SP,SPM} for some current works involving the geometry of operator spaces and orthogonality of operators in Banach space setting. Some recent developments on best approximations to compact operators can be found in \cite{Sai}, where the central themes are semi-inner-products \cite{G,L,Sa} and operator orthogonality.\\

The current article presents a comprehensive approach to address the problem of finding best approximation(s) to a given point $x$ out of a subspace $\mathbb{Y}$, in its full generality. After recalling some basic facts in Section \ref{Section:2}, we build up the theoretical background of our work in Section \ref{Section:3} and Section \ref{Section:4}. The results presented in Section \ref{Section:4} should be viewed as generalizations of the results obtained in \cite{Sai} and will be mentioned accordingly. The integral theme of our development is Birkhoff-James orthogonality of functionals. Application of Birkhoff-James orthogonality not only reduces the computational difficulties to resolve the above mentioned problem but also strengthens the classical duality principle \cite[Section 4]{Sai}. An extra advantage of employing the concept of Birkhoff-James orthogonality (over that of the classical duality principle) is that it provides an easy way out to compute the all possible best approximation(s) to $x$ out of $\mathbb{Y}$. Indeed, we devote Section \ref{Section:5} to show the applicability of the results, developed in the preceding sections, in context of the said problem. We obtain concrete solutions to some problems regarding best approximations and provide certain distance formulae under specific assumptions, which also produce some interesting inequalities, including a finite-dimensional strengthening of the classical H\"{o}lder's inequality. 

\section{Preliminaries}\label{Section:2}

In this section, we mention some known facts that will be used extensively in the next two sections. We begin with a simple proposition which has important applications in the study of topological vector spaces.

\begin{prop}\label{kernel and linear functionals}\cite[Lemma 3.9]{R}
Suppose that $g_1,g_2,\dots, g_m$ and $f$ are linear functionals on a vector space $\mathbb{X}.$ Let 
$$\mathcal{W}=\bigcap\limits_{i=1}^m \mathcal{N}(g_i)=\left\{x\in \mathbb{X}: g_1(x)=g_2(x)=\dots=g_m(x)=0\right\}.$$
Then the following three conditions are equivalent:\\

$(i)$ There exist scalars $\lambda_1, \lambda_2, \dots, \lambda_m$ such that
$$f=\lambda_1g_1+\lambda_2g_2+\dots +\lambda_mg_m.$$\\
$(ii)$ There exists $\gamma < \infty$ such that
$$|f(x)|\leq \gamma~ \underset{1\leq i\leq m}{\max}|g_i(x)|\qquad (x\in \mathbb{X}).$$\\
$(iii)$ $f(x)=0$ for every $x\in \mathcal{W}.$
\end{prop}

Suppose that $\tau$ is a topology on a vector space $\mathbb{X}$ such that every one point set in $\mathbb{X}$ is closed and the vector space operations on $\mathbb{X}$ are continuous with respect to the topology $\tau$. Then the vector space $\mathbb{X}$ equipped with the topology $\tau$ is called a topological vector space. The topological vector space $\mathbb{X}$ is called locally convex if there exists a local base at $\theta$, whose members are convex. 
Every topological vector space enjoys an important separation property:

\begin{prop}\label{Separation theorem}\cite[Theorem 1.10]{R}
Let $\mathbb{X}$ be a topological vector space. Let $K$ and $C$ be subsets of $\mathbb{X}$ such that $K$ is compact and $C$ is closed with $K\cap C= \emptyset$. Then there exists a neighborhood $V$ of $\theta$ such that
$$(K+V)\cap (C+V)=\emptyset$$ 
\end{prop}

Given a Banach space $\mathbb{X}$, $\mathbb{X}^*$ equipped with the weak$^*$ topology is a locally convex topological vector space. Moreover, every linear functional on $\mathbb{X}^*$ that is weak$^*$ continuous is of the form $\psi(x)$ for some $x\in \mathbb{X}$, where $\psi$ denotes the canonical embedding of $\mathbb{X}$ into $\mathbb{X}^{**}$. We refer the readers to the standard text \cite{R} for more information in this regard. Weak$^*$ topology on $\mathbb{X}^*$ has a crucial compactness property known as the Banach-Alaoglu Theorem:

\begin{theorem}(Banach-Alaoglu)\label{Banach-Alaoglu}
Let $\mathbb{X}$ be a normed linear space. Then the closed unit ball $B_{\mathbb{X}^*}$ of $\mathbb{X}^*$ is compact with respect to the weak$^*$ topology on $\mathbb{X}^*$.
\end{theorem}

We next present a classical result which is a variant of the geometric Hahn-Banach Theorem and is popularly known as the Mazur Theorem. In the following theorem, we do not require the topological vector space to be locally convex.

\begin{theorem}\label{Mazur}\cite[Theorem 18.2]{F}
Let $\mathbb{X}$ be a topological vector space and let $E$ be a linear subspace of $\mathbb{X}$. Let $V$ be a convex open subset of $\mathbb{X}$ such that 
$$E\cap V=\emptyset.$$ 
Then there exists a closed hyperplane $H$ of $\mathbb{X}$ such that 
$$E\subseteq H,\quad H\cap V=\emptyset.$$
\end{theorem}

Let $\mathbb{X}$ be a Banach space and let $\mathcal{W}$ be any non-trivial subspace of $\mathbb{X}$. Let $f$ be any member of $\mathbb{X}^*$. A member $f_0$ of $\mathbb{X}^*$ is said to be a \emph{Hahn-Banach extension} of $\restr{f}{\mathcal{W}}$ if 
$$\restr{f}{\mathcal{W}}=\restr{f_0}{\mathcal{W}}~\mathrm{and}~\left\|\restr{f}{\mathcal{W}}\right\|=\|f_0\|.$$\\

The next two propositions are about some basic facts regarding real $\ell_p$ spaces, where $p\in (1, \infty)$. Given any $1\leq p \leq \infty$, $q$ is said to be the conjugate to $p$ if $q=1$ $(q=\infty),$ whenever $p=\infty$ $(p=1),$ and, $q=\frac{p}{p-1},$ whenever $1<p<\infty.$

\begin{prop}\label{dual of lp is lq:1}
Let $p\in (1,\infty)$. Then the dual of $\ell_p$ is isometrically isomorphic to $\ell_q$, and for any $\mathbf{c}=(c_1,c_2,\dots)\in \ell_q$, the corresponding member $f_\mathbf{c}\in \ell_p^*$ is given by:
$$f_\mathbf{c}(x_1,x_2,\dots)=\sum\limits_{k\in \mathbb{N}}c_kx_k \qquad \forall~ (x_1,x_2,\dots)\in \ell_p.$$
\end{prop}

\begin{prop}\label{dual of lp is lq:2}
Let $p\in (1,\infty)$ and let $\mathbf{a}=(a_1,a_2, \dots)\in \ell_p$ be non-zero. Let $\mathbf{c}=(c_1,c_2,\dots)\in {\ell_q}$ be such that $\mathbf{c}$ corresponds to the support functional of $\mathbf{a}$ in $\ell_p^*$. Then for each $k\in \mathbb{N}$, $c_k$ is given by:
$$c_k=\frac{sgn(a_k)|a_k|^{p-1}}{\|\mathbf{a}\|_p^{p-1}}.$$
\end{prop}

\section{Orthogonality of functionals}\label{Section:3}

Our aim in this section is to obtain a characterization of Birkhoff-James orthogonality of weak$^*$ continuous functionals, which is of paramount importance in the context of our current work on best approximations. We would like to mention that the said characterization can also be obtained by modifying Theorem  2.1 of \cite{RSS}. However, we present a complete proof of the same, for the convenience of the readers. We need the following proposition, which completely describes the norm attainment sets of weak$^*$ continuous functionals, to serve our purpose.

\begin{prop}\label{Proposition 1 of norm attainment}
Let $\mathbb{X}$ be a Banach space and let $f\in \mathbb{X}^{**}$ be  weak$^*$ continuous. Then $M_f=\pm D$, where $D$ is a non-empty compact path connected subset of $B_{\mathbb{X}^*}$ with respect to the weak$^*$ topology on $\mathbb{X}^*$.
\end{prop}

\begin{proof}
Denote the canonical embedding of $\mathbb{X}$ into $\mathbb{X}^{**}$ by $\psi$. Since $f$ is weak$^*$ continuous, there exists $x_0\in \mathbb{X}$ such that $\psi(x_0)=f$. By the Hahn-Banach Theorem, there exists a support functional of $B_\mathbb{X}(\theta, \|x_0\|)$, say $x_0^*$, at $x_0$. Observe that
$$\psi(x_0)(x_0^*)=f(x_0^*)=x_0^*(x_0)=\|x_0\|=\|\psi(x_0)\|=\|f\|.$$
Therefore, $M_f\neq \emptyset$. Let $D$ be a subset of $B_{\mathbb{X}^*}$, defined by:
$$D=\{x^*\in B_{\mathbb{X^*}}: f(x^*)=\|f\|\}.$$
Evidently, $M_f=\pm D$ and $D$ is non-empty since $M_f$ is non-empty. Also, it is easy to see that $D$ is a closed subset of $B_{\mathbb{X}^*}$.
It now follows from the Banach-Alaoglu Theorem that $D$ is a compact subset of $B_{\mathbb{X}^*},$ with respect to the weak$^*$ topology on $\mathbb{X}^*$. Note that for any $x_1^*,x_2^*\in D$, $tx_1^*+(1-t)x_2^*\in D$ for all $t\in [0,1]$. Since $\mathbb{X}^*$ equipped with the weak$^*$ topology is a topological vector space, therefore, the mapping $t\mapsto tx_1^*+(1-t)x_2^*$ is continuous. Consequently, $D$ is a path connected subset of $B_{\mathbb{X}^*}$. This completes the proof.
\end{proof}

Based on Proposition \ref{Proposition 1 of norm attainment}, we now obtain a characterization of orthogonality of weak$^*$ continuous functionals. Note that the following characterization also includes \cite[Theorem 3.2]{Sai}.

\begin{theorem}\label{Orthogonality of functional}
Let $\mathbb{X}$ be a Banach space and let $f,g\in \mathbb{X}^{**}$ be weak$^*$ continuous. Then $f\perp_B g$ if and only if $M_f\cap \mathcal{N}(g)\neq \emptyset$.
\end{theorem}

\begin{proof}
We only prove the necessary part as the proof of the sufficient part is trivial. Suppose on the contrary that $M_f\cap \mathcal{N}(g)=\emptyset$. By Proposition \ref{Proposition 1 of norm attainment}, $M_f=\pm D$, where $D$ is a non-empty compact connected subset of $B_{\mathbb{X}^*}$ with respect to the weak$^*$ topology on $\mathbb{X}^*$. It now follows from the connectedness of $D$ that either $f(x)\cdot g(x)> 0$ for all $x\in D$, or, $f(x)\cdot g(x)< 0$ for all $x\in D$. Without loss of generality, let $f(x)\cdot g(x)> 0$ for all $x\in D$. Let $p:B_{\mathbb{X}^*}\times [-1,1]\to \mathbb{R}$ be defined by
$$p(x,\lambda)=|f(x)+\lambda g(x)|\qquad \forall~ (x, \lambda)\in B_{\mathbb{X}^*}\times [-1,1].$$
Obviously, $p$ is continuous. Also, for each $y\in M_f$ there exist a weak$^*$ open set $U_y$ containing $y$ and $\delta_y\in (0,1)$ such that
$$p(\tilde{y},\lambda)<\|f\|\qquad \forall~ (\tilde{y},\lambda)\in U_y\times (-\delta_y,0).$$
On the other hand, for each $z\in B_{\mathbb{X}^*}\setminus M_f$ there exist a weak$^*$ open set $V_z$ containing $z$ and $\delta_z\in (0,1)$ such that
$$p(\tilde{z},\lambda)<\|f\|\qquad \forall~ (\tilde{z},\lambda)\in V_z\times (-\delta_z,\delta_z).$$
Therefore, the collection $\{U_y:y\in M_f\}\cup \{V_z:z\in B_{\mathbb{X}^*}\setminus M_f \}$ forms a weak$^*$ open cover for $B_{\mathbb{X}^*}$. Due to the compactness of $B_{\mathbb{X}^*}$(with respect to the weak$^*$ topology on $\mathbb{X}^*$), there exist natural numbers $k_1,k_2$ such that 
$$B_{\mathbb{X}^*}\subseteq\left(\bigcup\limits_{i=1}^{k_1}U_{y_i}\right)\bigcup \left(\bigcup\limits_{j=1}^{k_2}V_{z_j}\right).$$
Let $0<\mu_0< \min\left\{\min\left\{\left\{\delta_{y_i}\right\}_{i=1}^{k_1}\right\}, \min\left\{\left\{\delta_{z_j}\right\}_{j=1}^{k_2}\right\}\right\}$. Since $f-\mu_0 g$ is weak$^*$ continuous, by Proposition \ref{Proposition 1 of norm attainment}, $M_{f-\mu_0g}\neq \emptyset$. Let $x_0\in M_{f-\mu_0g}$. Then it follows from the choice of $\mu_0$ that
$$\|f-\mu_0g\|=|(f-\mu_0g)(x_0)|<\|f\|.$$
This is a contradiction to the fact that $f\perp_B g$. This completes the proof.
\end{proof}

As an application of Theorem \ref{Orthogonality of functional} we have the following corollary:

\begin{cor}(James characterization of Birkhoff-James orthogonality)
Let $\mathbb{X}$ be a Banach space and let $x,y\in \mathbb{X}.$ Then $x\perp_B y$ if and only if there exists $x_0^*\in S_{\mathbb{X}^*}$ such that $x_0^*(x)=\|x\|$ and $x_0^*(y)=0$.
\end{cor}

\begin{proof}
Denote the canonical embedding of $\mathbb{X}$ into $\mathbb{X}^{**}$ by $\psi$. Let $\psi(x)=f$ and $\psi(y)=g$. Then $x\perp_B y$ if and only if $f\perp_B g$. Since $f,g\in \mathbb{X}^{**}$ are weak$^*$ continuous, it follows from Theorem \ref{Orthogonality of functional} that $x\perp_B y$ if and only if $M_f \cap \mathcal{N}(g)\neq \emptyset$. Let $x_0^*\in M_f \cap \mathcal{N}(g).$ Then
$$x_0^*(x)=f(x_0^*)=\|f\|=\|x\|~\mathrm{and}~x_0^*(y)=g(x_0^*)=0.$$
This completes the proof.
\end{proof}

\section{Birkhoff-James orthogonality and best approximations}\label{Section:4}

Let $\mathbb{X}$ be a reflexive, strictly convex Banach space. One of the fundamental ideas in \cite{Sai} was to identify the Banach space $\mathbb{X}$ with its double dual $\mathbb{X}^{**}$ and then treat the best approximation problem in $\mathbb{X}^{**}$ by employing Birkhoff-James orthogonality techniques. Unfortunately, the above idea does not work if $\mathbb{X}$ is not reflexive. This lacuna can be overcome by identifying $\mathbb{X}$ to the space of all weak$^*$ continuous functionals on $\mathbb{X}^*$. The following theorem provides a necessary and sufficient condition regarding the best approximation problem in the space of all weak$^*$ continuous functionals. In that sense, the result is a  generalization of Theorem 3.4 and Theorem 3.5 of \cite{Sai}. Also, note that we do not require the strict convexity of $\mathbb{X}$.

\begin{theorem}\label{Best approximation of functional}
Let $\mathbb{X}$ be a Banach space and let $f\in \mathbb{X}^{**}$ be weak$^*$ continuous. Let $\mathbb{Y}$ be a subspace of $\mathbb{X}^{**}$ such that each member of $\mathbb{Y}$ is weak$^*$ continuous and $f\notin \mathbb{Y}$. Let $g_0\in \mathbb{Y}$. Then $g_0$ is a best approximation to $f$ out of $\mathbb{Y}$ if and only if for every finite-dimensional subspace $\mathbb{Z}$ of $\mathbb{Y}$ containing $g_0$, $\bigcap\limits_{g\in \mathbb{Z}}\mathcal{N}(g)\bigcap M_{f-g_0} \neq \emptyset$.
\end{theorem}

\begin{proof}
We first prove the necessary part. Suppose on the contrary that there exists a  finite-dimensional subspace $\mathbb{Z}$ of $\mathbb{Y}$ containing $g_0$ such that 
$$\bigcap\limits_{g\in \mathbb{Z}}\mathcal{N}(g)\bigcap M_{f-g_0} =\emptyset.$$
Let $\left\{h_j:1\leq j \leq k\right\}$ be a basis of $\mathbb{Z}$. Consequently, $\bigcap\limits_{g\in \mathbb{Z}}\mathcal{N}(g)=\bigcap\limits_{j=1}^k\mathcal{N}(h_j).$
Also, it follows from Proposition \ref{Proposition 1 of norm attainment} that $M_{f-g_0}=\pm D$,  where $D$ is a compact convex subset of $B_{\mathbb{X}^*}$ with respect to the weak$^*$ topology on $\mathbb{X}^*.$ Since $\mathbb{X}^*$ equipped with the weak$^*$ topology is a locally convex  topological vector space and $\bigcap\limits_{g\in \mathbb{Z}}\mathcal{N}(g)$ is a closed subset of $\mathbb{X}^*$  disjoint from $D,$ there exists a convex neighbourhood $V$ of $\theta$ such that 
$$(D+V)\bigcap \left(\bigcap\limits_{g\in \mathbb{Z}} \mathcal{N}(g)+V\right)=\emptyset.$$ 
Since $D$ is convex, so is $D+V.$ In particular, $D+V$ is an open convex subset of $\mathbb{X}^*$, disjoint from $\bigcap\limits_{g\in \mathbb{Z}}\mathcal{N}(g).$ Therefore, it follows from Theorem \ref{Mazur} that there exists a closed hyperplane $H$ of $\mathbb{X}^*$ such that 
$$\bigcap\limits_{g\in \mathbb{Z}}\mathcal{N}(g)\subseteq H,\quad H\cap (D+V)=\emptyset.$$ 
Let $h:\mathbb{X}^* \to \mathbb{R}$ be a linear functional such that $\mathcal{N}(h)=H.$ It now follows from Proposition \ref{kernel and linear functionals} that $h\in span\left\{h_j:1\leq j \leq k\right\}\subseteq \mathbb{Y}.$ In particular, $h$ is weak$^*$ continuous. Also, note that 
$$M_{f-g_0} \cap \mathcal{N}(h) =\emptyset.$$
Theorem \ref{Orthogonality of functional} ensures that $f-g_0\not\perp_B h.$ However, this is a contradiction to the fact that $g_0$ is a best approximation to $f$ out of $\mathbb{Y}.$\\

To prove the sufficient part of the theorem, let $g_1\in \mathbb{Y}$ be arbitrary and let $\mathbb{Z}=span\{g_0,g_1\}.$ Then it follows from the hypothesis of the theorem that 
$$\bigcap\limits_{g\in \mathbb{Z}}\mathcal{N}(g)\bigcap M_{f-g_0}\neq \emptyset.$$
In particular, $M_{f-g_0}\cap \mathcal{N}(g_1)\neq \emptyset.$ It now follows from Theorem \ref{Orthogonality of functional} that $f-g_0\perp_B g_1$. Since $g_1\in \mathbb{Y}$ was chosen arbitrarily, we have that $(f-g_0)\perp_B \mathbb{Y}.$ Consequently, $g_0\in \mathbb{Y}$ is a best approximation to $f$ out of $\mathbb{Y}$ and the proof follows.
\end{proof}

Whenever $\mathbb{Y}$ is finite-dimensional, the above theorem takes a simpler form. We record this as a corollary. The proof of the corollary follows directly from Theorem \ref{Best approximation of functional}, and therefore, it is omitted.

\begin{cor}\label{Best approximation of functional: corollary}
Let $\mathbb{X}$ be a Banach space and let $f,g_1,g_2, \dots, g_m\in \mathbb{X}^{**}$ be weak$^*$ continuous. Let $g_1,g_2, \dots, g_m$ be linearly independent and $f\notin \mathbb{Y}$, where $\mathbb{Y}=span\{g_1,g_2, \dots, g_m\}$. Then $g_0\in \mathbb{Y}$ is a best approximation to $f$ out of $\mathbb{Y}$ if and only if $\bigcap\limits_{i=1}^m\mathcal{N}(g_i)\bigcap M_{f-g_0} \neq \emptyset$.
\end{cor}

As an immediate application of Theorem \ref{Best approximation of functional}, we now obtain a distance formula in the space of all weak$^*$ continuous functionals on $\mathbb{X}^*$. The following distance formula can be regarded as a strengthened version of Theorem 3.6 of \cite{Sai}.

\begin{theorem}\label{Distance formula}
Let $\mathbb{X}$ be a Banach space and let $f\in \mathbb{X}^{**}$ be weak$^*$ continuous. Let $\mathbb{Y}$ be a subspace of $\mathbb{X}^{**}$ such that each member of $\mathbb{Y}$ is weak$^*$ continuous and $f\notin \mathbb{Y}$. Suppose that $g_0\in \mathbb{Y}$ is a best approximation to $f$ out of $\mathbb{Y}$. Then for any finite-dimensional subspace $\mathbb{Z}$ of $\mathbb{Y}$ containing $g_0$,
\begin{align*}
\|f-g_0\|=dist(f,\mathbb{Y})=\max \left\{ |f(x)|: x\in \bigcap\limits_{g\in \mathbb{Z}}\mathcal{N}(g)\bigcap S_{\mathbb{X}^*}\right\}.
\end{align*}
\end{theorem}

\begin{proof}
It follows from Theorem \ref{Best approximation of functional} that $\bigcap\limits_{g\in \mathbb{Z}}\mathcal{N}(g)\bigcap M_{f-g_0}\neq \emptyset$. Also, note that 
\begin{align}\label{restriction equation}
\restr{\left(f-g_0\right)}{\bigcap\limits_{g\in \mathbb{Z}}\mathcal{N}(g)}=f.
\end{align}
Fix some $x_0\in \bigcap\limits_{g\in \mathbb{Z}}\mathcal{N}(g)\bigcap M_{f-g_0}$. Then it is easy to see that
$$\|f-g_0\|=|(f-g_0)(x_0)|=|f(x_0)|\qquad \mathrm{(using~(\ref{restriction equation}))}.$$
We now claim that 
$$|f(x_0)|=\max \left\{ |f(x)|: x\in \bigcap\limits_{g\in \mathbb{Z}}\mathcal{N}(g)\bigcap S_{\mathbb{X}^*}\right\}.$$
Indeed, if there exists $y_0\in \bigcap\limits_{g\in \mathbb{Z}}\mathcal{N}(g)\bigcap S_{\mathbb{X}^*}$ such that
$|f(y_0)|>|f(x_0)|$, then we obtain that
$$\|f-g_0\|\geq |f(y_0)|>|f(x_0)| =|(f-g_0)(x_0)|=\|f-g_0\|,$$
which is a contradiction. This completes the proof. 
\end{proof}

Assuming $\mathbb{Y}$ to be finite-dimensional in the above theorem, we have the following distance formula:

\begin{cor}\label{Distance formula: corollary}
Let $\mathbb{X}$ be a Banach space and let $f,g_1,g_2, \dots, g_m\in \mathbb{X}^{**}$ be weak$^*$ continuous for some $m\in \mathbb{N}$. Let $g_1,g_2, \dots, g_m$ be linearly independent and let $f\notin \mathbb{Y}$, where $\mathbb{Y}=span\{g_1,g_2, \dots, g_m\}$. Then
\begin{align*}
dist(f,\mathbb{Y})=\max \left\{ |f(x)|: x\in \bigcap\limits_{i=1}^m\mathcal{N}(g_i)\bigcap S_{\mathbb{X}^*}\right\}.
\end{align*}
\end{cor}

\begin{proof}
A standard compactness argument ensures that there exist $\alpha_i\in \mathbb{R}$, where $1\leq i \leq m$, such that $\sum\limits_{i=1}^m \alpha_ig_i$ is a best approximation to $f$ out of $\mathbb{Y}$. Now, arguing as in Theorem \ref{Distance formula}, we get the desired formula.
\end{proof}

As mentioned in the introduction, Birkhoff-James orthogonality techniques provide some genuine insights in determining best approximation(s) to a given point out of a subspace. Indeed, the following result completely characterizes best approximation(s) to a given point out of a finite-dimensional subspace in the setting of weak$^*$ continuous functionals. 

\begin{theorem}\label{Finding best approximation}
Let $\mathbb{X}$ be a Banach space and let $f,f_0,g_1,g_2, \dots, g_m\in \mathbb{X}^{**}$ be weak$^*$ continuous for some $m\in \mathbb{N}$. Let $g_1,g_2, \dots, g_m$ be linearly independent and let $f\notin \mathbb{Y}$, where $\mathbb{Y}=span\{g_1,g_2, \dots, g_m\}$. Then $f-f_0$ is a best approximation to $f$ out of $\mathbb{Y}$ if and only if $f_0$ is a Hahn-Banach extension of $\restr{f}{\mathcal{W}}$, where $\mathcal{W}=\bigcap\limits_{i=1}^m \mathcal{N}(g_i)$.
\end{theorem}

\begin{proof}
We first prove the necessary part. Suppose on the contrary that $\restr{f}{\mathcal{W}}\neq \restr{f_0}{\mathcal{W}}.$ Then there exists $x_0\in \mathcal{W}$ such that $f(x_0)\neq f_0(x_0)$. However, then it follows from Proposition \ref{kernel and linear functionals} that $f-f_0\notin \mathbb{Y}$, which is a contradiction. Next, suppose that $\left\|\restr{f}{\mathcal{W}}\right\|\neq \|f_0\|.$ Then $\left\|\restr{f}{\mathcal{W}}\right\|< \|f_0\|.$ We now claim that $M_{f_0}\cap \mathcal{W}=\emptyset$. Indeed, if $y_0\in M_{f_0}\cap \mathcal{W}$, then it follows that
$$\|f_0\|=|f_0(y_0)|=|f(y_0)|\leq \left\|\restr{f}{\mathcal{W}}\right\|<\|f_0\|,$$
which is absurd. Therefore, $M_{f_0}\cap \mathcal{W}=\emptyset$, as expected. Then Corollary \ref{Best approximation of functional: corollary} ensures that $f-f_0$ is not a best approximation to $f$ out of $\mathbb{Y}$, which is a contradiction.\\ 

We now prove the sufficient part. Since $\restr{f}{\mathcal{W}}= \restr{f_0}{\mathcal{W}}$, $\mathcal{N}(f-f_0)$ contains $\mathcal{W}$. It now follows from Proposition \ref{kernel and linear functionals} that $f-f_0\in \mathbb{Y}$. Next, consider any $y_1\in M_{\restr{f}{\mathcal{W}}}$. Then
$$\left\|\restr{f}{\mathcal{W}}\right\|=\left|f(y_1)\right|=|f_0(y_1)|=\|f_0\|.$$
Thus, $y_1\in M_{f_0}$ and $M_{\restr{f}{\mathcal{W}}}\subseteq M_{f_0}$. Consequently, $M_{f_0}\cap \mathcal{W}\neq \emptyset.$ Therefore, applying Corollary \ref{Best approximation of functional: corollary}, we obtain that $f-f_0$ is a best approximation to $f$ out of $\mathbb{Y}$. This completes the proof.
\end{proof}

We would like to remark here that the above theorem can also be stated in terms of norm attainment sets of weak$^*$ continuous functionals.

\begin{remark}
Let $\mathbb{X}$ be a Banach space and let $f,g_1,g_2, \dots, g_m\in \mathbb{X}^{**}$ be weak$^*$ continuous for some $m\in \mathbb{N}$. Let $g_1,g_2, \dots, g_m$ be linearly independent and let $f\notin \mathbb{Y}$, where $\mathbb{Y}=span\{g_1,g_2, \dots, g_m\}$. Let us consider the set 
$$\mathcal{B}:=\{h\in \mathbb{X}^{**}: f-h \mathrm{~is~a~best~approximation~to}~f\mathrm{~out~of}~\mathbb{Y}\}.$$
Then $\mathcal{B}$ is precisely the collection of those weak$^*$ continuous functionals which are the extensions of $\restr{f}{\bigcap\limits_{i=1}^m \mathcal{N}(g_i)}$ and whose norm attainment sets contain the norm attainment set of the restriction of $f$ to $\bigcap\limits_{i=1}^m \mathcal{N}(g_i).$ In other words, for any $h\in \mathcal{B}$,
$$\restr{h}{\mathcal{W}}=\restr{f}{\mathcal{W}}~\mathrm{and}~M_{\restr{f}{\mathcal{W}}}\subseteq M_{h},~\mathrm{where}~ \mathcal{W}=\bigcap\limits_{i=1}^m \mathcal{N}(g_i).$$
\end{remark}

It is obvious that given any Banach space $\mathbb{X}$, the best approximation problems in $\mathbb{X}$ can be treated as the best approximation problems in the space of all weak$^*$ continuous functionals on $\mathbb{X}^*$. Therefore, the foregoing results can also be rephrased in terms of best approximation problems in $\mathbb{X}$. As an evidence of this, we present the following theorem which is essentially a variant of Theorem \ref{Best approximation of functional}. However, we refrain ourselves from doing analogous treatment to the remaining results, to avoid monotony.

\begin{theorem}\label{Best approximation in domain space}
Let $\mathbb{X}$ be a Banach space and let $x_0\in \mathbb{X}$. Let $\mathbb{Y}$ be a subspace of $\mathbb{X}$ such that $x_0\notin \mathbb{Y}$ and $y_0\in \mathbb{Y}$. Then $y_0$ is a best approximation to $x_0$ out of $\mathbb{Y}$ if and only if for every finite-dimensional subspace $\mathbb{Z}$ of $\mathbb{Y}$ containing $y_0$, $\bigcap\limits_{z\in \mathbb{Z}}\mathcal{N}(\psi(z))\bigcap M_{\psi(x_0-y_0)} \neq \emptyset$, where $\psi:\mathbb{X}\to \mathbb{X}^{**}$ denotes the canonical embedding.
\end{theorem}

We have seen in Theorem \ref{Orthogonality of functional} that there is a deep connection  between Birkhoff-James orthogonality and the norm attainment set of a given functional. Since every weak$^*$ continuous functional is norm attaining, there is a scope to employ Birkhoff-James orthogonality techniques in all of the preceding results. However, Birkhoff-James orthogonality is not so straightforward for the functionals that do not attain their norms. Consequently, the same techniques cannot be applied for functionals not attaining their respective norms. We end this section with a result which addresses this issue. Also, note that the result is valid in any Banach space, \emph{real or complex}.

\begin{theorem}\label{Best approximations for infinite-dimensional subspaces}
Let $\mathbb{X}$ be a Banach space and let $f\in \mathbb{X}^*$. Let $\mathbb{Y}$ be a subspace of $\mathbb{X}^{*}$ such that $f\notin \mathbb{Y}$ and $g_0\in \mathbb{Y}.$ Then $g_0$ is a best approximation to $f$ out of $\mathbb{Y}$ if and only if for every finite-dimensional subspace $\mathbb{Z}$ of $\mathbb{Y}$ containing $g_0$, the following holds true:
\begin{align*}
\left\|\restr{(f-g_0)}{\bigcap\limits_{g\in \mathbb{Z}}\mathcal{N}(g)}\right\|=\|f-g_0\|.
\end{align*}
\end{theorem}

\begin{proof}
We first prove the necessary part. Suppose on the contrary that there exists a finite-dimensional subspace $\mathbb{Z}$ of $\mathbb{Y}$ containing $g_0$ such that
\begin{align*}
\left\|\restr{(f-g_0)}{\bigcap\limits_{g\in \mathbb{Z}}\mathcal{N}(g)}\right\|<\|f-g_0\|.
\end{align*}
Let $\left\{h_j:1\leq j \leq k\right\}$ be a basis of $\mathbb{Z}$. Then it is straightforward to check that $\bigcap\limits_{g\in \mathbb{Z}}\mathcal{N}(g)=\bigcap\limits_{j=1}^k\mathcal{N}(h_j).$  
Let $f_0:\bigcap\limits_{j=1}^k \mathcal{N}\left(h_{j}\right)\to \mathbb{R}$ be defined by
$$f_0(x)=f(x) \qquad \forall~ x\in \bigcap\limits_{j=1}^k \mathcal{N}\left(h_{j}\right).$$
By the Hahn-Banach Theorem, $f_0$ possesses a linear extension $\widetilde{f_0}:\mathbb{X}\to \mathbb{R}$ such that
$$\|\widetilde{f_0}\|=\left\|{f_0}\right\|=\left\|\restr{f}{\bigcap\limits_{j=1}^k\mathcal{N}(h_{j})}\right\|.$$
Now, consider the linear functional $(f-\widetilde{f_0}):\mathbb{X}\to \mathbb{R}$. Since $(f-\widetilde{f_0})$ vanishes identically on $\bigcap\limits_{j=1}^k\mathcal{N}(h_{j})$, it follows from Proposition \ref{kernel and linear functionals} that $(\widetilde{f_0}-f)\in span\{h_{j} :1\leq j \leq k\}\subseteq \mathbb{Y}.$ On the other hand,
$$\|f-(f-\widetilde{f_0})\|=\|\widetilde{f_0}\|=\left\|\restr{f}{\bigcap\limits_{j=1}^k\mathcal{N}(h_{j})}\right\|=\left\|\restr{f}{\bigcap\limits_{g\in \mathbb{Z}}\mathcal{N}(g)}\right\|=\left\|\restr{(f-g_0)}{\bigcap\limits_{g\in \mathbb{Z}}\mathcal{N}(g)}\right\|<\|f-g_0\|,$$
where the second last equality follows from the fact that $f(x)=(f-g_0)(x)$ for all $x\in \bigcap\limits_{g\in \mathbb{Z}}\mathcal{N}(g)$. However, this is a contradiction to the fact that $g_0$ is a best approximation to $f$ out of $\mathbb{Y}$.\\

We now prove the sufficient part. Let $g_1\in \mathbb{Y}$ be arbitrary and let $\mathbb{Z}=span\{g_0,g_1\}$. Then it follows from the hypothesis of the theorem that
\begin{align*}
\|f-g_1\|\geq \left\|\restr{(f-g_1)}{\bigcap\limits_{g\in \mathbb{Z}}\mathcal{N}(g)}\right\|=\left\|\restr{(f-g_0)}{\bigcap\limits_{g\in \mathbb{Z}}\mathcal{N}(g)}\right\|=\|f-g_0\|.
\end{align*}
Thus, $g_0$ is a best approximation to $f$ out of $\mathbb{Y}$ and this completes the proof.
\end{proof}

\section{Some applications and distance formulae}\label{Section:5}

In this section we exhibit some interesting applications and examples to the theories developed in the preceding sections. Let us begin with an algorithm that generalizes the Problem given in \cite{Sai}.

\begin{problem}\label{Problem:1}
Let $\mathbb{X}$ be a Banach space and let $m\in \mathbb{N}$. Let $\mathbf{x_0},\mathbf{y}_1,\mathbf{y}_2, \dots, \mathbf{y}_m\in \mathbb{X}$ be such that $\mathbf{y}_1,\mathbf{y}_2, \dots, \mathbf{y}_m$ are linearly independent and $\mathbf{x_0}\notin \mathbb{Y}$, where $\mathbb{Y}=span\{\mathbf{y}_1,\mathbf{y}_2, \dots, \mathbf{y}_m\}$. Then find best approximation(s) to $\mathbf{x_0}$ out of $\mathbb{Y}$ and compute $dist(\mathbf{x_0},\mathbb{Y})$.
\end{problem}

Corollary \ref{Best approximation of functional: corollary}, Corollary \ref{Distance formula: corollary} and Theorem \ref{Finding best approximation} allow us to approach the problem in the following three steps:\\

\noindent Step 1: We embed $\mathbb{X}$ into its double dual $\mathbb{X}^{**}$ via the canonical isometric isomorphism $\psi$. Let $\psi(\mathbf{x_0})=f$ and $\psi(\mathbf{y}_i)=g_i$ for $1\leq i\leq m$. Let $\mathbb{Z}= span\{g_1,g_2, \dots, g_m\}$. Evidently, the above problem is equivalent to finding the best approximation(s) to $f$ out of $\mathbb{Z}$ and computing $dist(f,\mathbb{Z})$. Since $f,g_1,g_2,\dots,g_m$ are weak$^*$ continuous, the criteria of Corollary \ref{Best approximation of functional: corollary}, Corollary \ref{Distance formula: corollary} and Theorem \ref{Finding best approximation} are satisfied.\\

\noindent Step 2: Let $\mathcal{W}=\bigcap\limits_{i=1}^m \mathcal{N}(g_i)$. We now consider the following two cases:\\
\noindent Case I: $\mathcal{W}$ is one-dimensional. Consider any non-zero $z\in \mathcal{W}$. Then it follows from Corollary \ref{Distance formula: corollary} that 
\begin{align*}
dist(\mathbf{x_0}, \mathbb{Y})=dist(f,\mathbb{Z})=\max \left\{ |f(x)|: x\in \mathcal{W}\cap S_{\mathbb{X}^*}\right\}=\frac{1}{\|z\|}|f(z)|.\\
\end{align*}

\noindent Case II: $\mathcal{W}$ is not one-dimensional. Consider $\mathcal{W}\cap kerf$ and find some non-zero $u\in \mathcal{W}$ such that $u\perp_B \left(\mathcal{W}\cap kerf\right)$. Note that the existence of such an $u$ is always guaranteed, since $M_{\restr{f}{\mathcal{W}}}\neq \emptyset$. Then we have that
$$\max \left\{ |f(x)|: x\in \mathcal{W}\cap S_{\mathbb{X}^*}\right\}=\frac{1}{\|u\|}|f(u)|.$$
In other words,
\begin{align*}
dist(\mathbf{x_0}, \mathbb{Y})=dist(f,\mathbb{Z})=\max \left\{ |f(x)|: x\in \mathcal{W}\cap S_{\mathbb{X}^*}\right\}=\frac{1}{\|u\|}|f(u)|.\\
\end{align*}

\noindent Step 3: Let us consider the following collection:
$$\Lambda:=\left\{f-f_0\in \mathbb{X}^{**}:~f_0~\mbox{is~a~Hahn-Banach~extension~of}~\restr{f}{\mathcal{W}}\right\}.$$
It is not difficult to see that $\Lambda \neq \emptyset$ and $\Lambda \subseteq \mathbb{Z}.$ It follows from Theorem \ref{Finding best approximation} that $\Lambda$ is precisely the collection of best approximation(s) to $f$ out of $\mathbb{Z}$.\\

Thus, we completely obtain the solution of the above problem, as $\psi^{-1}(\Lambda)$ is precisely the collection of best approximation(s) to $\mathbf{x_0}$ out of $\mathbb{Y}$.\\

The algorithm presented in Problem \ref{Problem:1} is particularly advantageous for $\ell_1^n$ and $\ell_\infty^n$ spaces, for $n\in \mathbb{N}.$ This is because the dual of $\ell_1^n$ ($\ell_\infty^n$) is $\ell_\infty^n$  ($\ell_1^n$) and if any member $x^*$ of $\ell_1^{n^*}$ ($\ell_\infty^{n^*}$) corresponds to a member $(a_1,a_2,\dots, a_n)$ of $\ell_\infty^n$ ($\ell_1^n$), then the action of $x^*$ on any member $(b_1,b_2, \dots, b_n)$ of $\ell_1^n$ ($\ell_\infty^n$) is given by the formula:
$$x^*(b_1,b_2,\dots, b_n)=\sum\limits_{i=1}^n a_ib_i.$$

We elaborate this in more detail in the following problem:

\begin{problem}\label{Problem:2}
Let $\mathbb{X}=\ell_1^n$ or $\ell_\infty^n$ for some $n\in \mathbb{N}$. Let $\mathbf{x_0},\mathbf{y}_1,\mathbf{y}_2,\dots,\mathbf{y}_m\in \mathbb{X}$, where $1\leq m <n$ be such that $\mathbf{y}_1,\mathbf{y}_2,\dots,\mathbf{y}_m$ are linearly independent and $\mathbf{x_0}\notin \mathbb{Y}$, where $\mathbb{Y}=span\{\mathbf{y}_1,\mathbf{y}_2,\dots,\mathbf{y}_m\}$. Compute $dist(\mathbf{x_0}, \mathbb{Y})$.
\end{problem}

\noindent Step 1: Let $\mathbf{x_0}=(x_1,x_2,\dots, x_n)$ and let $\mathbf{y}_i=(y_{i1},y_{i2},\dots,y_{in})$, where $1\leq i \leq m$. Let $\Omega(\mathbf{x_0})=f$ and $\Omega(\mathbf{y}_i)=g_i$, where $1\leq i \leq m$. Let $\mathbb{Z}=span\{g_1,g_2,\dots,g_m\}$. Evidently, the above problem is equivalent to computing $dist(f,\mathbb{Z})$.\\

\noindent Step 2: Note that
$$\bigcap\limits_{i=1}^m\mathcal{N}(g_i)=\mathcal{W}=\left\{(u_1,u_2, \dots, u_n)\in \mathbb{R}^n: \sum\limits_{j=1}^n y_{ij}u_j = 0; 1\leq i \leq m\right\}.$$ 
Therefore, applying Corollary \ref{Distance formula: corollary}, we obtain
$$dist(\mathbf{x_0},\mathbb{Y})=dist(f,\mathbb{Z})=\max \left\{ |f(x)|:x\in \mathcal{W}\bigcap S_{\mathbb{X}^*}\right\}.$$
In other words, the problem of computing $dist(\mathbf{x_0}, \mathbb{Y})$ reduces to the problem of finding the absolute maximum of $f$ on the unit sphere of the solution space of the system of linear equations: 
$$\sum\limits_{j=1}^n y_{ij}u_j=0;\qquad 1\leq i \leq m.$$
Note that the non-triviality of the solution space is guaranteed by the existence of best approximation(s).\\

We now present an example to illustrate the utility of the above problem:

\begin{example}\label{Example to algorithm}
Let $\mathbb{X}=\ell_1^4$ and let $\mathbb{Y}=span\{\mathbf{y}_1,\mathbf{y}_2\}\subseteq \mathbb{X}$, where $\mathbf{y}_1=(1,2,0,0)$, $\mathbf{y}_2=(-1,0,2,0).$ Let $\mathbf{x_0}=(1,1,1,1).$ Then it is trivial to see that $\mathbf{x_0}\notin \mathbb{Y}.$ Our aim is to calculate $dist(\mathbf{x_0}, \mathbb{Y})$ and to find a best approximation to $\mathbf{x_0}$ out of $\mathbb{Y}$.\\

Let $\Omega: \ell_1^4 \to \ell_\infty^{4^*}$ denote the canonical isometric isomorphism. Let 
$$\Omega(\mathbf{x_0})=f~\mathrm{and}~\Omega(\mathbf{y}_1)=g_1,\Omega(\mathbf{y}_2)=g_2.$$
Let $\mathbb{Z}=span\{g_1,g_2\}.$ It is easy to see that
$$\mathcal{N}(g_1)\cap \mathcal{N}(g_2)=span\left\{(1,-\frac{1}{2},\frac{1}{2},0),(0,0,0,1)\right\}.$$
We now follow the procedure, as described in Case II of Problem \ref{Problem:1}. Therefore, we find $\bigcap\limits_{i=1}^2 \mathcal{N}(g_i)\bigcap \mathcal{N}(f),$ which is given by
$span\left\{(1, -\frac{1}{2}, \frac{1}{2},-1)\right\}.$ Note that $(1, -\frac{1}{2}, \frac{1}{2},1)\in \bigcap\limits_{i=1}^2 \mathcal{N}(g_i)\bigcap S_{\ell_\infty^4}$ and $(1, -\frac{1}{2}, \frac{1}{2},1)\perp_B (1, -\frac{1}{2}, \frac{1}{2},-1).$ Thus,
$$|f(1, -\frac{1}{2}, \frac{1}{2},1)|=dist(\mathbf{x_0},\mathbb{Y})=dist(f,\mathbb{Z})=\max \left\{ |f(x)|:x\in \bigcap\limits_{i=1}^2 \mathcal{N}(g_i)\bigcap S_{\mathbb{X}^*}\right\}=2.$$
Next, we find a best approximation to $f$ out of $\mathbb{Z}.$ Define $f_0: \bigcap\limits_{i=1}^2 \mathcal{N}(g_i)\to \mathbb{R}$ by:
$$f_0(x_1,x_2,x_3,x_4)=\restr{f}{\bigcap\limits_{i=1}^2 \mathcal{N}(g_i)}(x_1,x_2,x_3,x_4)=x_1+x_4\qquad \forall~(x_1,x_2,x_3,x_4)\in \bigcap\limits_{i=1}^2 \mathcal{N}(g_i).$$
Let $\widetilde{f_0}:\ell_\infty^4\to \mathbb{R}$ be defined by
$$\widetilde{f_0}(x_1,x_2,x_3,x_4)=x_1+x_4 \qquad \forall~(x_1,x_2,x_3,x_4)\in \ell_\infty^4.$$
Then we have that
$$\|\widetilde{f_0}\|=\left\|{f_0}\right\|=\left\|\restr{f}{\bigcap\limits_{i=1}^2\mathcal{N}(g_{i})}\right\|=2.$$
Consequently, $\widetilde{f_0}$ is a Hahn-Banach extension of $\restr{f}{\bigcap\limits_{i=1}^2\mathcal{N}(g_{i})}.$
Therefore, it follows from Theorem \ref{Finding best approximation} that $f-\widetilde{f_0}$ is a best approximation to $f$ out of $\mathbb{Z}.$ Thus, $\Omega^{-1}(f-\widetilde{f_0})=(0,1,1,0)$ is a best approximation to $\mathbf{x_0}$ out of $\mathbb{Y}.$
\end{example}

In light of Problem \ref{Problem:2}, we can say that best approximation problems in context of the $\ell_\infty^n$ and $\ell_1^n$ spaces reduce to the problem of maximizing a functional to the unit sphere of the solution space of a system of homogeneous linear equations, and the  problem becomes trivial if $m=n-1$. Indeed, in that case all we need to do is to solve a system of homogeneous linear equations. On the other hand, we have the following explicit distance formulae whenever $n=2$:

\begin{prop}\label{Distance of l_1}
Let $ \mathbb{X} = \ell_1^2, $ and let $ \mathbf{x} = (a, b) \in \mathbb{X} $. Let $\mathbf{y}=(c,d)\neq (0,0)$ be such that $\mathbf{x}\notin \mathbb{Y}$, where $\mathbb{Y} = span\{ \mathbf{y} \}.$ Then
\begin{align*}
dist (\mathbf{x}, \mathbb{Y}) = \frac{2|ad-bc|}{{|d|+|c|}+{\left||d|-|c|\right|}}.
\end{align*}
\end{prop}

\begin{proof}
Let $ \Omega : \ell_1^2 \longrightarrow \ell_{\infty}^{{2}^{*}} $ be the canonical isometric isomorphism. Let $ \Omega(\mathbf{x}) = f $ and let $ \Omega(\mathbf{y}) = g. $ Clearly,
$$\mathcal{N}(g) = \left\{(u_1, u_2) \in \mathbb{R}^2 : cu_1+du_2=0 \right\}.$$
Assume that $d \neq 0.$ Then it follows that $ (1, -\frac{c}{d}) \in \mathcal{N}(g). $ Let $ \alpha \in \mathbb{R} $ be such that $\left\| \alpha\left (1, -\frac{c}{d}\right) \right\|_{\infty} = 1.$ It can be shown without any difficulty that 
$$|\alpha| = \dfrac{2|d|}{{|d|+|c|}+{||d|-|c||}}.$$ Applying Corollary \ref{Distance formula: corollary}, we obtain that 
$$dist(\mathbf{x}, \mathbb{Y}) = \left| f\left(\alpha \left(1, -\frac{c}{d}\right)\right) \right| = \frac{2|ad-bc|}{{|d|+|c|}+{\left||d|-|c|\right|}}.$$
Similarly, if $ c \neq 0, $ one can show that 
$$dist (\mathbf{x}, \mathbb{Y}) = \frac{2|ad-bc|}{{|d|+|c|}+{\left||d|-|c|\right|}}.$$
This completes the proof.
\end{proof}

\begin{prop}\label{Distance of l_infinity}
Let $ \mathbb{X} = \ell_\infty^2, $ and let $ \mathbf{x} = (a, b) \in \mathbb{X} $. Let $\mathbf{y}=(c,d)\neq (0,0)$ be such that $\mathbf{x}\notin \mathbb{Y}$, where $\mathbb{Y} = span\{ \mathbf{y} \}.$ Then
\begin{align*}
dist(\mathbf{x}, \mathbb{Y}) = \frac{|ad-bc|}{|c|+|d|}.
\end{align*}
\end{prop}

\begin{proof}
Let $ \Omega : \ell_\infty^2 \longrightarrow \ell_{1}^{{2}^{*}} $ be the canonical isometric isomorphism. Let $ \Omega(\mathbf{x}) = f $ and let $ \Omega(\mathbf{y}) = g. $ Clearly,
$$\mathcal{N}(g) = \left\{(u_1, u_2) \in \mathbb{R}^2 : cu_1+du_2=0 \right\}.$$
Assume that $d \neq 0.$ Then it follows that $ (1, -\frac{c}{d}) \in \mathcal{N}(g). $ Let $ \alpha \in \mathbb{R} $ be such that $\left\| \alpha\left (1, -\frac{c}{d}\right) \right\|_1 = 1.$ It can be shown without any difficulty that $|\alpha| = \dfrac{|d|}{|c|+|d|}.$ Applying Corollary \ref{Distance formula: corollary}, we obtain that 
$$dist(\mathbf{x}, \mathbb{Y}) = \left| f\left(\alpha \left(1, -\frac{c}{d}\right)\right) \right| = \frac{|ad-bc|}{|c|+|d|}.$$
Similarly, if $ c \neq 0, $ one can show that 
$$dist (\mathbf{x}, \mathbb{Y}) = \frac{|ad-bc|}{|c|+|d|}.$$
This completes the proof.
\end{proof}

In view of Proposition \ref{Distance of l_1}, Proposition \ref{Distance of l_infinity} and Theorem $3.8$ of \cite{Sai},  the proof of the following result is obvious.

\begin{prop}\label{Distance of l_1, l_infty, l_p}
Let $ \mathbb{X} = \ell_p^2,$ $1\leq p \leq \infty$; and let $ \mathbf{x} = (a, b) \in \mathbb{X} $. Let $\mathbf{y}=(c,d)\neq (0,0)$ be such that $\mathbf{x}\notin \mathbb{Y}$, where $\mathbb{Y} = span\{ \mathbf{y} \}.$ Then
\begin{align*}
dist (\mathbf{x}, \mathbb{Y}) = \frac{|ad-bc|}{\|(c,d)\|_q},
\end{align*}
where $q$ is conjugate to $p$.
\end{prop}

The next theorem presents a sufficient condition for the uniqueness of the best approximation in finite-dimensional real polyhedral Banach space. Recall that a finite-dimensional real Banach space is called \emph{polyhedral} if $Ext(B_\mathbb{X})$ is finite. For more information on the local structure of finite-dimensional real polyhedral Banach spaces, we refer the readers to \cite{SPBB}.

\begin{theorem}\label{Best approximation in polygonal}
Let $\mathbb{X}$ be an $n$-dimensional real polyhedral Banach space and let $x,y_1,y_2,\dots, y_m\in \mathbb{X};$ $1\leq m <n$. Let $y_1,y_2,\dots, y_m$ be linearly independent and $x\notin \mathbb{Y}$, where $\mathbb{Y}=span\{y_1,y_2,\dots,y_m\}$. Let $\psi:\mathbb{X}\to \mathbb{X}^{**}$ denote the canonical isometric isomorphism. Suppose that $\bigcap\limits_{i=1}^m\mathcal{N}(\psi(y_i)) \bigcap S_{\mathbb{X}^*}$ contains only smooth point(s) of $B_{\mathbb{X}^*}$. Then the best approximation to $x$ out of $\mathbb{Y}$ is unique.
\end{theorem}

\begin{proof}
Let $\psi(x)=f$ and let $\psi(y_i)=g_i$ for each $1\leq i \leq m$. Clearly, $f\notin \mathbb{Z}$, where $\mathbb{Z}= span\{g_1,g_2, \dots, g_m\}$. It is enough to show that the best approximation to $f$ out of $\mathbb{Z}$ is unique. Suppose that $\sum\limits_{i=1}^m \alpha_ig_i$ and $\sum\limits_{i=1}^m \beta_ig_i$ are best approximations to $f$ out of $\mathbb{Z}$ for real numbers $\alpha_1, \alpha_2, \dots, \alpha_m;\beta_1, \beta_2,\dots, \beta_m$. Note that $\bigcap\limits_{i=1}^m\mathcal{N}(g_i) \bigcap S_{\mathbb{X}^*}\subseteq F\cup (-F)$ for some facet $F$ of $B_{\mathbb{X}^*}$. Indeed, if $x_1\in  \left(\bigcap\limits_{i=1}^m\mathcal{N}(g_i) \bigcap S_{\mathbb{X}^*}\right)\bigcap F_1$ and $x_2\in  \left(\bigcap\limits_{i=1}^m\mathcal{N}(g_i) \bigcap S_{\mathbb{X}^*}\right)\bigcap F_2$ for some distinct facets $F_1$ and $F_2$ of $B_{\mathbb{X}^*}$ with $F_1\neq -F_2$, then there exists $t\in (0,1)$ such that $\frac{tx_1+(1-t)x_2}{\|tx_1+(1-t)x_2\|}\in \bigcap\limits_{i=1}^m\mathcal{N}(g_i) \bigcap S_{\mathbb{X}^*}$ is a non-smooth point of $B_{\mathbb{X}^*}$. Let $g_0$ be a best approximation to $f$ out of $\mathbb{Z}$. Then it follows from Corollary \ref{Best approximation of functional: corollary} that
$$\bigcap\limits_{i=1}^m\mathcal{N}(g_i) \bigcap M_{f-g_0}\neq \emptyset.$$
Since $\bigcap\limits_{i=1}^m\mathcal{N}(g_i) \bigcap S_{\mathbb{X}^*}$ contains only smooth points of the facets $F$ and $-F$, it follows that $f-g_0=\lambda h$, where $h\in S_{\mathbb{X}^{**}}$ is the unique support functional corresponding to the facet $F$ \cite{SPBB} and $\lambda\in \mathbb{R}$ is non-zero. However, this proves that
$$f-\sum\limits_{i=1}^m \alpha_ig_i=\lambda_1 h~\mathrm{and}~f-\sum\limits_{i=1}^m \beta_ig_i= \lambda_2 h,$$
where $\lambda_1,\lambda_2$ are non-zero real numbers with $|\lambda_1|=|\lambda_2|$. Suppose that $\lambda_1=-\lambda_2=\lambda_0$. Then it follows that
$$f=\frac{1}{2}\sum\limits_{i=1}^m(\alpha_i+ \beta_i)g_i,$$
which is a contradiction, since $f\notin \mathbb{Z}$. Therefore, we must have $\lambda_1=\lambda_2$. Consequently, $\sum\limits_{i=1}^m \beta_ig_i=\sum\limits_{i=1}^m \alpha_ig_i$ and the best approximation is unique. This completes the proof of the theorem.
\end{proof}

The converse of the above theorem need not be true. The following example illustrates such a situation:

\begin{example}\label{Non-strictly convex example}
Let $\mathbb{X}=\ell_1^3$ and let $x=(0,\frac{1}{2}, \frac{1}{2})$, $y=(0,0,1)$. We identify $\mathbb{X}^{**}$ to the dual of $\ell_\infty^3$ and let $\psi:\mathbb{X}\to \mathbb{X}^{**}$ denote the canonical isometric isomorphism. Let $\psi(x)=f$ and $\psi(y)=g$. Then it is not difficult to see that
$$f(x_1,x_2,x_3)=\frac{1}{2}(x_2+x_3)~\mathrm{and}~g(x_1,x_2,x_3)=x_3 \qquad \forall~ (x_1,x_2,x_3)\in \ell_\infty^3.$$
Clearly,
$$\mathcal{N}(g)\cap S_{\mathbb{X}^*}=\left\{(x_1,x_2,0):x_1,x_2\in \mathbb{R}, \max\{|x_1|,|x_2|\}=1\right\},$$
which also contains non-smooth points of $B_{\mathbb{X}^*}.$ Our aim is to show that the best approximation to $x$ out of $span\{y\}$ is unique.\\

Let $f-\lambda_0g$ be the best approximation to $f$ out of $span\{g\}$. Then it follows from Corollary \ref{Distance formula: corollary} that
\begin{align*}
dist\left(f,span\{g\}\right)=\max \left\{ \frac{1}{2}|(x_2+x_3)|: x_3=0,\max\{|x_1|,|x_2|\}=1\right\}=\frac{1}{2}.
\end{align*}
Define $f_0: \mathcal{N}(g)\to \mathbb{R}$ by
$$f_0(x_1,x_2,x_3)=\restr{f}{\mathcal{N}(g)}(x_1,x_2,x_3)=\frac{1}{2}x_2\qquad \forall~ (x_1,x_2,x_3)\in \mathcal{N}(g).$$
Next, we define $\widetilde{f_0}:\ell_\infty^3\to \mathbb{R}$ by
$$\widetilde{f_0}(x_1,x_2,x_3)=\frac{1}{2}x_2\qquad \forall~ (x_1,x_2,x_3)\in \ell_\infty^3.$$
Now, it is not difficult to see that
$$\|\widetilde{f_0}\|=\|f_0\|=\left\|\restr{f}{\mathcal{N}(g)}\right\|=\frac{1}{2}.$$
Consequently, $\widetilde{f_0}$ is a Hahn-Banach extension of $\restr{f}{\mathcal{N}(g)}.$ Therefore, it follows from Theorem \ref{Finding best approximation} that $f-\widetilde{f_0}=\frac{1}{2}g$ is a best approximation to $f$ out of $span\{g\}.$\\

For any $\lambda\in \mathbb{R}$, consider the linear functional $(f-\lambda g):\ell_\infty^3\to \mathbb{R}$, given by:
$$(f-\lambda g)(x_1,x_2,x_3)=\frac{1}{2}x_2+(\frac{1}{2}-\lambda)x_3 \qquad \forall~ (x_1,x_2,x_3)\in \ell_\infty^3.$$
Observe that if $(\frac{1}{2}-\lambda)>0$, choosing $(x_0,y_0,z_0)\in Ext(B_{\mathbb{X}^*})$ such that $y_0z_0=1$, we get $|(f-\lambda g)(x_0,y_0,z_0)|>\frac{1}{2}$. Again, if $(\frac{1}{2}-\lambda)<0$, choosing $(x_0,y_0,z_0)\in Ext(B_{\mathbb{X}^*})$ such that $y_0z_0=-1$, we get $|(f-\lambda g)(x_0,y_0,z_0)|>\frac{1}{2}$. Therefore, the only solution of $\lambda_0$ for which $\|f-\lambda_0g\|=\frac{1}{2}$ is $\lambda_0=\frac{1}{2}$ and the best approximation to $f$ out of $span\{g\}$ is unique. Consequently, the best approximation to $x$ out of $span\{y\}$ is unique.
\end{example}


As an application of the ideas developed in Theorem \ref{Best approximation of functional} and Theorem \ref{Distance formula: corollary}, it is also possible to explore the following invariant distance problem:

\begin{problem}
Let $(\mathbb{X}_1,\|\cdot\|_1)$ and $(\mathbb{X}_2,\|\cdot\|_2)$ be two Banach spaces with $\mathbb{X}_1\subseteq \mathbb{X}_2$. Let $\mathbb{Y}$ be a finite-dimensional vector subspace of $\mathbb{X}_1\cap \mathbb{X}_2$ and let $x_0\in (\mathbb{X}_1\cap \mathbb{X}_2)\setminus \mathbb{Y}$. Then find a necessary and sufficient condition on $x_0$ and $\mathbb{Y}$ such that
$$dist_1(x_0, \mathbb{Y})=dist_2(x_0, \mathbb{Y}).$$
\end{problem}

In the following theorem, we provide a complete solution to the above problem when $\mathbb{X}_1=\ell_{p_1}$ and $\mathbb{X}_2=\ell_{p_2}$, where $1<p_1, p_2 <\infty$ and $p_1\neq p_2.$

\begin{theorem}\label{equality of distance}
Let $1<p_1<p_2<\infty$. Let $\mathbf{x}, \mathbf{y}_1,\dots, \mathbf{y}_m\in \ell_{p_1}\cap \ell_{p_2}$ be such that $\mathbf{y}_1,\mathbf{y}_2,\dots, \mathbf{y}_m$ are linearly independent with $\mathbf{x}\notin \mathbb{Y},$ where $\mathbb{Y}=span\{\mathbf{y}_1,\mathbf{y}_2,\dots, \mathbf{y}_m\}.$ Let $\mathbf{x_0}$ and $\mathbf{y_0}$ be the best approximations to $\mathbf{x}$ out of $\mathbb{Y}$ in $\ell_{p_1}$ and $\ell_{p_2}$, respectively. Then $dist_{p_1}(\mathbf{x},\mathbb{Y})=dist_{p_2}(\mathbf{x}, \mathbb{Y})$ if and only if $\mathbf{x}-\mathbf{x_0}=\mathbf{x}-\mathbf{y_0}=\lambda e_j$, for some non-zero $\lambda\in \mathbb{R}$ and $j\in \mathbb{N}$, and $\mathbb{Y}\subseteq \left\{(\lambda_1, \lambda_2, \dots)\in \ell_{p_1}\cap \ell_{p_2}: \lambda_j=0\right\}.$ 
\end{theorem}

\begin{proof}
We only prove the necessary part as the sufficient part of the theorem is trivial. Let $\Omega_1:\ell_{p_1}\to \ell_{q_1}^*$ and $\Omega_2:\ell_{p_2}\to \ell_{q_2}^*$ denote the canonical isometric isomorphisms, where $q_1$ and $q_2$ are conjugates to $p_1$ and $p_2$, respectively. Given any $\eta\in \ell_{p_1}$, $\eta$ is also a member of $\ell_{p_2}$ and by Proposition \ref{dual of lp is lq:1}, 
\begin{align}\label{equality of functionals in duals}
\restr{\Omega_1(\eta)}{\ell_{q_2}}=\Omega_2(\eta).
\end{align}
Let $\mathbb{Z}_1=span\{\Omega_1(\mathbf{y}_i): 1\leq i \leq m\}$ and $\mathbb{Z}_2=span\{\Omega_2(\mathbf{y}_i): 1\leq i \leq m\}.$ It follows from the hypothesis of the theorem that $dist_{q_1^*}(\Omega_1(\mathbf{x}), \mathbb{Z}_1)=dist_{q_2^*}(\Omega_2(\mathbf{x}), \mathbb{Z}_2).$ Let
$$\bigcap\limits_{i=1}^m\mathcal{N}(\Omega_1(\mathbf{y}_i))=\mathcal{W}_1, \quad \bigcap\limits_{i=1}^m\mathcal{N}(\Omega_2(\mathbf{y}_i))=\mathcal{W}_2.$$
Applying Corollary \ref{Distance formula: corollary}, we have that
\begin{align*}
dist_{q_1^*}(\Omega_1(\mathbf{x}), \mathbb{Z}_1)& =\max \left\{ |\Omega_1(\mathbf{x})(z)|: z\in \mathcal{W}_1\cap S_{\ell_{q_1}}\right\}\\
& =dist_{q_2^*}(\Omega_2(\mathbf{x}), \mathbb{Z}_2)\\
&=\max \left\{ |\Omega_2(\mathbf{x})(z)|: z\in \mathcal{W}_2\cap S_{\ell_{q_2}}\right\}\\
& = \lambda~\mathrm{(say)}.    
\end{align*}
Since $\mathcal{W}_1$ and $\mathcal{W}_2$ are strictly convex, there exist $u_0\in S_{\ell_{q_1}}$ and $v_0\in S_{\ell_{q_2}}$ such that
$$M_{\restr{\Omega_1(\mathbf{x})}{\mathcal{W}_1}}=\{\pm u_0\}, \quad M_{\restr{\Omega_2(\mathbf{x})}{\mathcal{W}_2}}=\{\pm v_0\}.$$
Obviously, $|\Omega_1(\mathbf{x})(u_0)|=|\Omega_2(\mathbf{x})(v_0)|.$
Also, it is not difficult to see that $v_0\in \ell_{q_1}$ and $\|v_0\|_{q_1}\leq 1.$ Now, it follows from (\ref{equality of functionals in duals}) that $\Omega_1(\mathbf{y}_i)(v_0)=\restr{\Omega_1(\mathbf{y}_i)}{\ell_{q_2}}(v_0)=\Omega_2(\mathbf{y}_i)(v_0)=0$ for each $1\leq i \leq m.$ Therefore, $v_0\in \mathcal{W}_1.$ If $\|v_0\|_{q_1}< 1,$ consider some $\mu_0>1$ such that $\mu_0v_0\in S_{\ell_{q_1}}$. Then
$$|\Omega_1(\mathbf{x})(\mu_0v_0)|=\mu_0|\Omega_1(\mathbf{x})(v_0)|=\mu_0|\Omega_2(\mathbf{x})(v_0)|=\mu_0|\Omega_1(\mathbf{x})(u_0)|>|\Omega_1(\mathbf{x})(u_0)|,$$
which is a contradiction. Therefore, $\|v_0\|_{q_1}=1.$ Also, $u_0=\pm v_0$, as otherwise $\restr{\Omega_1(\mathbf{x})}{\mathcal{W}_1}$ would attain norm at two pair of points. Without loss of generality, let $u_0=v_0$ and $\Omega_1(\mathbf{x})(u_0)=\Omega_2(\mathbf{x})(v_0)=\lambda.$ Since $q_1\neq q_2$, we must have $u_0=v_0=\pm e_j$ for some $j\in \mathbb{N}.$ Let $u_0=v_0=e_j$ and let $h_1$, $h_2$ be the (unique) best approximations to $\Omega_1(\mathbf{x})$, $\Omega_2(\mathbf{x})$ out of $\mathbb{Z}_1$ and $\mathbb{Z}_2$, respectively. Obviously, $\Omega_1(\mathbf{x}_0)=h_1$ and $\Omega_2(\mathbf{y}_0)=h_2.$ By Theorem \ref{Finding best approximation}, $\Omega_1(\mathbf{x}-\mathbf{x}_0)$ and $\Omega_2(\mathbf{x}-\mathbf{y}_0)$ are the Hahn-Banach extensions of $\restr{\Omega_1(\mathbf{x})}{\mathcal{W}_1}$ and $\restr{\Omega_2(\mathbf{x})}{\mathcal{W}_2}$, respectively. Since $\ell_{q_1}$ and $\ell_{q_2}$ are strictly convex, we have $M_{\Omega_1(\mathbf{x}-\mathbf{x}_0)}=M_{\Omega_2(\mathbf{x}-\mathbf{y}_0)}=\{\pm e_j\}.$ Therefore, $\frac{1}{\lambda}\Omega_1(\mathbf{x}-\mathbf{x}_0)$ and $\frac{1}{\lambda}\Omega_2(\mathbf{x}-\mathbf{y}_0)$ are the (unique) support functionals at $e_j$. Again it follows from Proposition \ref{dual of lp is lq:2} that
\begin{align*}
&\frac{1}{\lambda}\Omega_1(\mathbf{x}-\mathbf{x}_0)(z)=z_j~\quad \forall~z=(z_1, z_2, \dots)\in \ell_{q_1},\\
&\frac{1}{\lambda}\Omega_2(\mathbf{x}-\mathbf{y}_0)(z)=z_j~\quad \forall~z=(z_1, z_2, \dots)\in \ell_{q_2}.
\end{align*}
Now, applying Proposition \ref{dual of lp is lq:1}, we have that
$$\mathbf{x}-\mathbf{x_0}=\mathbf{x}-\mathbf{y_0}=\lambda e_j.$$ 
Let $\rho= J\{e_j\}$. Since $e_j\perp_B \mathbb{Y}$, we have
$$\rho(y)=y_j=0 \quad \forall~y=(y_1,y_2, \dots)\in \mathbb{Y}.$$ Consequently, $\mathbb{Y}\subseteq \left\{(\lambda_1, \lambda_2, \dots)\in \ell_{p_1}\cap \ell_{p_2}: \lambda_j=0\right\}$ and the proof follows.
\end{proof}

It can be seen from \cite[Remark 3.10]{Sai} that best approximation problems give rise to a family of inequalities in context of $\ell_p$ spaces $(1<p<\infty)$. Our next goal is to find the said family of inequalities in a more general setting. The following result is the first step towards achieving the said goal.

\begin{theorem}\label{Inequalities}
Let $n\in \mathbb{N}$ and let $\mathbb{X}_i=\ell_{p_i}^{m_i}$, where $m_i$ are natural numbers and $1\leq p_i \leq \infty$ for each $1\leq i \leq n$. Let
$$\mathbb{X}=\bigoplus\limits_{i=1}^n{_p~\mathbb{X}_i}\qquad \textit{for~some}~1\leq p \leq \infty.$$ 
Let $\sum\limits_{i=1}^n m_i=m$ and let $T$ be an $m$ by $m$ non-zero matrix. Then for any $\mathbf{x_0}\in \mathbb{X}$ the following holds true:
\begin{align}\label{Optimal inequalities}
\min\{\|\mathbf{x_0}-y\|_\mathbb{X}: y\in \mathcal{R}(T)\}=\max\left\{ \frac{1}{\|z\|_{\mathbb{X}^*}}\left|\langle \mathbf{x_0}, z\rangle\right|: z\in \mathcal{N}(T^*)\setminus \{\theta\}\right\},
\end{align}
\end{theorem}

\begin{proof}
We begin the proof with the observation that the existence of the minimum of the set $\{\|\mathbf{x_0}-y\|_\mathbb{X}: y\in \mathcal{R}(T)\}$ is guaranteed by the existence of best approximation(s) to $\mathbf{x_0}$ out of $\mathcal{R}(T)$.\\
Let $q,q_1,q_2, \dots, q_n$ be the conjugates to $p,p_1,p_2, \dots,p_n$, respectively. Clearly, for each $1\leq i \leq n$,
$$\mathbb{X}_i^*=\ell_{q_i}^{m_i},~\mathrm{and}~\mathbb{X}^*=\bigoplus\limits_{i=1}^n{_q~\mathbb{X}_i^*}.$$
Also, for any member $x^*=(x_1^*, x_2^*, \dots, x_n^*)\in \mathbb{X}^*$,
$$x^*(x_1, x_2, \dots x_n)=\sum\limits_{i=1}^n x_i^*x_i \qquad \forall~ (x_1, x_2, \dots x_n)\in \bigoplus\limits_{i=1}^n{_p~\mathbb{X}_i}.$$
Clearly, $\mathbb{X}$ and $\mathbb{R}^m$ are isomorphic as vector spaces. Let $\{e_j\}_{j=1}^m$ denote the standard ordered basis of $\mathbb{X}$. Observe that $\mathcal{R}(T)=span\{\mathbf{y}_1,\mathbf{y}_2,\dots, \mathbf{y}_m\}$, where $\mathbf{y}_j=T(e_j)$ for all $1\leq j \leq m$. Therefore,
$$[T]=[\mathbf{y}_1^t~ \mathbf{y}_2^t~ \dots ~\mathbf{y}_m^t].$$
Let $\mathbf{x_0}=(u_1,u_2, \dots u_n)$, where $u_i\in \mathbb{X}_i$ for all $1\leq i \leq n$. Also, let
$$\mathbf{y}_j=(w_{j1},w_{j2},\dots, w_{jn}),\quad \mathrm{where}\quad w_{ji}\in\mathbb{X}_i \qquad \forall~ 1\leq i \leq n.$$
Let $\Omega: \bigoplus\limits_{i=1}^n{_p~\mathbb{X}_i}\to \left(\bigoplus\limits_{i=1}^n{_q~\mathbb{X}_i^*}\right)^*$ denote the canonical isometric isomorphism. Let
$$\Omega(\mathbf{x_0})=f_0\quad \mathrm{and}\quad \Omega(\mathbf{y}_j)=g_j \qquad \forall ~1\leq j \leq m.$$
Then it is not difficult to see that for each $1\leq j \leq m$,
\begin{align*}
& g_j(x_1,x_2, \dots, x_n)= \sum\limits_{i=1}^n \langle w_{ji}, x_i \rangle \qquad & \forall ~(x_1,x_2, \dots, x_n)\in \bigoplus\limits_{i=1}^n{_q~\mathbb{X}_i^*},
\end{align*}
\begin{align}\label{lp functional}
& f_0(x_1,x_2, \dots, x_n)= \sum\limits_{i=1}^n \langle u_i, x_i \rangle \qquad & \forall ~(x_1,x_2, \dots, x_n)\in \bigoplus\limits_{i=1}^n{_q~\mathbb{X}_i^*}.
\end{align}
Let $\mathbb{Z}=span\{g_1, g_2, \dots, g_m\}$. Note that $\bigcap\limits_{j=1}^m \mathcal{N}(g_j)$ is a subspace of $\bigoplus\limits_{i=1}^n \mathbb{R}^{m_i}$. Identifying each member $z=(z_1,z_2, \dots, z_n)$ of $\bigoplus\limits_{i=1}^n{_q\mathbb{X}_i^*}$ as a member of $\bigoplus\limits_{i=1}^n \mathbb{R}^{m_i}\simeq\mathbb{R}^m$, we then have
$$\bigcap\limits_{j=1}^m \mathcal{N}(g_j)=\left\{(z_1,z_2, \dots, z_n)\in \bigoplus\limits_{i=1}^n \mathbb{R}^{m_i}:\sum\limits_{i=1}^n \langle w_{1i}, z_i \rangle =\dots=\sum\limits_{i=1}^n \langle w_{mi}, z_i \rangle =0\right\}.$$
In other words,
$$\bigcap\limits_{j=1}^m \mathcal{N}(g_j)= \left\{z\in \mathbb{R}^m: [\mathbf{y}_1^t~ \mathbf{y}_2^t~ \dots ~\mathbf{y}_m^t]^tz= \theta\right\}= \left\{z\in \mathbb{R}^m : T^*z=\theta\right\}=\mathcal{N}(T^*).$$

We now consider the following two cases:\\

Case I: Let $\mathbf{x_0}\in \mathcal{R}(T).$ Then $\min\{\|\mathbf{x_0}-y\|_\mathbb{X}: y\in \mathcal{R}(T)\}=0$. Since $\mathbf{x_0}\in span\{\mathbf{y}_1, \mathbf{y}_2, \dots, \mathbf{y}_m\}$, we have that $f_0\in \mathbb{Z}$. It now follows from Proposition \ref{kernel and linear functionals} that
$$\mathcal{N}(T^*)=\bigcap\limits_{j=1}^m \mathcal{N}(g_j)\subseteq kerf_0.$$
Therefore, $\max \left\{|f_0(z)|: z\in \bigcap\limits_{j=1}^m \mathcal{N}(g_j)\bigcap S_{\bigoplus\limits_{i=1}^n{_q\mathbb{X}_i^*}}\right\}=0.$ Now, applying (\ref{lp functional}), we have that $\max \left\{\frac{1}{\|z\|_{\mathbb{X}^*}}|\langle \mathbf{x_0}, z\rangle| : z\in \mathcal{N}(T^*)\setminus \{\theta\}\right\}=0$ and the equality (\ref{Optimal inequalities}) follows.\\

\noindent Case II: Let $\mathbf{x_0}\notin \mathcal{R}(T)$ and let $y\in \mathcal{R}(T)$. Then applying Corollary \ref{Distance formula: corollary} we have
\begin{align*}
\|\mathbf{x_0}-y\|_\mathbb{X}  = \left\|\Omega(\mathbf{x_0}-y)\right\|_{\mathbb{X}^{**}}  &= \left\|f_0-\Omega(y)\right\|_{\mathbb{X}^{**}} \geq dist(f_0, \mathbb{Z})\\
& = \max \left\{|f_0(z)|: z\in \bigcap\limits_{j=1}^m \mathcal{N}(g_j)\bigcap S_{\bigoplus\limits_{i=1}^n{_q\mathbb{X}_i^*}}\right\}\\ 
& = \max \left\{\frac{1}{\|z\|_{\mathbb{X}^*}}|f_0(z)|: z\in \mathcal{N}(T^*)\setminus \{\theta\}\right\}.
\end{align*}
Now, applying (\ref{lp functional}), we have that
\begin{align*}
\|\mathbf{x_0}-y\|_\mathbb{X} \geq \max \left\{\frac{1}{\|z\|_{\mathbb{X}^*}}|\langle \mathbf{x_0}, z\rangle| : z\in \mathcal{N}(T^*)\setminus \{\theta\}\right\}. 
\end{align*}
Since $y\in \mathcal{R}(T)$ was chosen arbitrarily, we obtain
$$\min\{\|\mathbf{x_0}-y\|_\mathbb{X}: y\in \mathcal{R}(T)\}\geq \max \left\{\frac{1}{\|z\|_{\mathbb{X}^*}}|\langle \mathbf{x_0}, z\rangle| : z\in \mathcal{N}(T^*)\setminus \{\theta\}\right\}.$$
However, the above inequality is necessarily an equality, since
$$dist(\mathbf{x_0},\mathcal{R}(T))=dist(f_0, \mathbb{Z})=\min\{\|\mathbf{x_0}-y\|_\mathbb{X}: y\in \mathcal{R}(T)\}.$$
This completes the proof.
\end{proof}

The above result can also be stated in the form of the following inequality:

\begin{cor}\label{Application of Inequalities}
Let $n\in \mathbb{N}$ and let $\mathbb{X}_i=\ell_{p_i}^{m_i}$, where $m_i$ are natural numbers and $1\leq p_i \leq \infty$ for each $1\leq i \leq n$. Let
$$\mathbb{X}=\bigoplus\limits_{i=1}^n{_p~\mathbb{X}_i}\qquad \textit{for~some}~1\leq p \leq \infty.$$ 
Let $\sum\limits_{i=1}^n m_i=m$ and let $T$ be an $m$ by $m$ non-zero matrix. Let $\mathbf{x_0}\in \mathbb{X}$. Then for any $y\in \mathcal{R}(T)$
\begin{align*}
\|\mathbf{x_0}-y\|_\mathbb{X}\geq \frac{1}{\|z\|_{\mathbb{X}^*}}|\langle \mathbf{x_0}, z\rangle|,
\end{align*}
for all $z\in \mathcal{N}(T^*)\setminus \{\theta\}.$
\end{cor}

As an application of Theorem \ref{Inequalities}, we have the following:

\begin{theorem}\label{Hilbert space and inequalities}
Let $\mathbf{x_0}=(x_1,x_2, \dots, x_m)\in \mathbb{R}^m$ and let $\{m_1, m_2, \dots ,m_n\}\subseteq \mathbb{N}$ be such that $m=\sum\limits_{i=1}^n m_i$. Set $s_0=0$ and $s_k=\sum\limits_{i=1}^k m_i$ for each $1\leq k \leq n$. Let $\left\{p,p_1, \dots, p_n\right\}\subseteq (1, \infty)$ and let $q,q_1,q_2, \dots, q_n$ are conjugates to $p,p_1,p_2, \dots, p_n,$ respectively. Then for any $\mathbf{a}=(a_1,a_2, \dots, a_m)\in \mathbb{R}^m\setminus \{\theta\}$
\begin{align}\label{Hilbert space inequalities: expression}
&\left(\mathlarger{\mathlarger{\sum}}_{k=1}^n\left(\sum\limits_{j=s_{k-1}+1}^{s_k}|x_j-\lambda a_j|^{p_k}\right)^{\dfrac{p}{p_k}}\right)^{\dfrac{1}{p}}\left(\mathlarger{\mathlarger{\sum}}_{k=1}^n\left(\sum\limits_{j=s_{k-1}+1}^{s_k}|b_j|^{q_k}\right)^{\dfrac{q}{q_k}}\right)^{\dfrac{1}{q}}\geq \left | \sum\limits_{j=1}^m x_jb_j\right|,
\end{align}
for all $(b_1,b_2, \dots, b_m)\in W\setminus \{\theta\}$ and $\lambda \in \mathbb{R}$, where
$$W=\left\{(z_1,z_2, \dots, z_m)\in \mathbb{R}^m: \sum\limits_{i=1}^m a_iz_i=0\right\}.$$
Moreover, the above inequality is optimal.
\end{theorem}

\begin{proof}
Let $\mathbb{H}$ denote the Hilbert space $\mathbb{R}^m$ equipped with the usual dot product. Let $\mathbb{Y}=span\{(a_1,a_2, \dots, a_m)\}$. Therefore, we have that $W=\mathbb{Y}^{\perp}$. Let for each $1\leq i \leq n$
$$\mathbb{X}_i=\ell_{p_i}^{m_i}~\mathrm{and}~\mathbb{X}=\bigoplus\limits_{i=1}^n{_p~\mathbb{X}_i}.$$
 Then $\mathbb{X}$ and $\mathbb{R}^m$ are isomorphic as vector spaces. Let $\{e_j\}_{j=1}^m$ denote the standard ordered basis of $\mathbb{X}$. Define $T:\mathbb{H}\to \mathbb{H}$ in such a way that
$$\mathcal{R}(T)=span\{T(e_j): 1\leq j \leq m\}=\mathbb{Y}.$$
Then it is easy to see that $\mathcal{N}(T^*)=\mathbb{Y}^\perp$. It now follows from Theorem \ref{Inequalities} that
\begin{align*}
\min\{\|\mathbf{x_0}-y\|_\mathbb{X}: y\in \mathcal{R}(T)\}=\max\left\{ \frac{1}{\|z\|_{\mathbb{X}^*}}\left|\langle \mathbf{x_0}, z\rangle\right|: z\in \mathcal{N}(T^*)\setminus \{\theta\}\right\},
\end{align*}
Now, inequality (\ref{Hilbert space inequalities: expression}) is obtained by removing minimum and maximum from both sides, and expressing $\|\mathbf{x_0}-y\|_\mathbb{X}$ and $\|z\|_{\mathbb{X}^*}$, in their explicit forms. Moreover, the inequality is optimal. This completes the proof.
\end{proof}

The inequality (\ref{Hilbert space inequalities: expression}) obtained in Theorem \ref{Hilbert space and inequalities} exhibits a stronger version of H\"{o}lder's inequality in finite-dimensional case. We explain this in more detail in the following remark:
\begin{remark}\label{Holder inequality}
We consider two non-zero elements $u=(u_1,u_2, \dots,u_m)$ and $v=(v_1,v_2,\dots, v_m)$ in $\mathbb{R}^m$, with $m\geq 2$.
Let $\mathbf{x_0}=(u_1,u_2, \dots, u_m)$ and $\widetilde{v}=(b_1,b_2,\dots,b_m)$, where
$$b_j=\begin{cases}
sgn({u_j}{v_j})v_j &~\mathrm{if}~u_j\neq 0,~1\leq j \leq m,\\
v_j &  ~\mathrm{if}~u_j= 0,~1\leq j \leq m.
\end{cases}$$
Let $\mathbf{a}=(a_1,a_2, \dots, a_m)\in \mathbb{R}^m\setminus \{\theta\}$ be such that $\sum\limits_{j=1}^m a_jb_j=0$. Let $\{m_1, m_2, \dots ,m_n\}$ be any subset of natural numbers such that $m=\sum\limits_{i=1}^n m_i$. Now, considering $p=p_1=p_2=\dots=p_n$ and $\lambda=0$ in Theorem \ref{Hilbert space and inequalities}, we obtain that
\begin{align*}
\left(\sum\limits_{j=1}^m|u_j|^p\right)^{\dfrac{1}{p}}\left(\sum\limits_{j=1}^m|b_j|^q\right)^{\dfrac{1}{q}}\geq \left | \sum\limits_{j=1}^m u_jb_j\right|.
\end{align*}
Since $\sum\limits_{j=1}^m|b_j|^q=\sum\limits_{j=1}^m|v_j|^q$ and $\sum\limits_{j=1}^m u_jb_j=\sum\limits_{j=1}^m |u_jv_j|$, on simplification, we get
\begin{align*}
\sum\limits_{j=1}^m |u_jv_j| \leq \left(\sum\limits_{j=1}^m|u_j|^p\right)^{\dfrac{1}{p}}\left(\sum\limits_{j=1}^m|v_j|^q\right)^{\dfrac{1}{q}}.
\end{align*}
\end{remark}

We end this section with an example involving a particular type of minimization problem. It is worth mentioning that such kind of problems are difficult to handle from the algebraic point of view. However, employing the duality techniques developed in this work from the perspective of orthogonality, we can immediately solve these problems via trivial computations.

\begin{example}\label{Example: Minimization problem}
Let $\mathbf{x}_0=(\alpha_1, \alpha_2, \alpha_3, \alpha_4, \alpha_5, \alpha_6, \alpha_7, \alpha_8, \alpha_9, \alpha_{10})\in \mathbb{R}^{10}$ and let $\{a_i\}_{i=1}^{10}$ be real variables. Let us consider the following minimization problem:
\begin{align}\label{Minimization problem}
& \min_{\substack {A_i\\ 1\leq i \leq 10}} \left\{\begin{array}{l} |\alpha_1-A_1|^5+ \left(|\alpha_2-A_2|^7+|\alpha_3-A_3|^7\right)^{\dfrac{5}{7}}\\
+\left(|\alpha_4-A_4|^3+|\alpha_5-A_5|^3+|\alpha_6-A_6|^3\right)^{\dfrac{5}{3}}\\
+\left(|\alpha_7-A_7|^{11}+|\alpha_8-A_8|^{11}\right)^{\dfrac{5}{11}}+\left(|\alpha_9-A_9|^9+|\alpha_{10}-A_{10}|^9\right)^{\dfrac{5}{9}}\end{array}\right\}^{\dfrac{1}{5}},
\end{align}

where $A_i$ ($1\leq i \leq 10)$ are given by:
\begin{align*}
& A_1 = -9a_1+7a_2+9a_3+9a_4+8a_5+9a_6+6a_7+7a_8-8a_9+8a_{10},\\
& A_2 = a_1-4a_3+a_4+3a_5+3a_6-2a_7+6a_8+9a_{10},\\
& A_3 = a_1+5a_2+3a_3+6a_4+4a_5-8a_6+9a_7+7a_8+2a_{10},\\
& A_4 = 3a_1+9a_2-4a_3+3a_4+8a_5+2a_6+5a_7+9a_8+a_9+7a_{10},\\
& A_5 = 6a_1+5a_2+7a_3-a_4+6a_5+a_6+8a_7+8a_8-6a_9+5a_{10},\\
& A_6 = 8a_1+3a_2+8a_3-a_4+2a_6+a_7+3a_8+9a_9+5a_{10},\\
& A_7 = 2a_2+8a_3-7a_4+3a_5+4a_7+2a_8+6a_{10},\\
& A_8 = 7a_1+7a_2+a_3+5a_5+2a_6+a_7+a_8+4a_9+9a_{10},\\
& A_9 = 9a_1+a_2+9a_3+5a_4-3a_5+7a_6+5a_7+3a_8+9a_9+8a_{10},\\
& A_{10} = 12a_1+\frac{13}{2}a_2+13a_3-\frac{3}{2}a_4+\frac{5}{2}a_5+\frac{11}{2}a_6+\frac{11}{2}a_7+\frac{9}{2}a_8+11a_9+14a_{10}.
\end{align*}

In view of Theorem \ref{Inequalities}, we obtain the solution of the above problem in the following three steps:\\

\noindent Step I: The above minimization problem has \emph{five summands}:
\begin{align*}
& \mathrm{first~summand:} |\alpha_1-A_1|^5\\
& \mathrm{second~summand:} \left(|\alpha_2-A_2|^7+|\alpha_3-A_3|^7\right)^{\dfrac{5}{7}},\\
& \mathrm{third~summand:} \left(|\alpha_4-A_4|^3+|\alpha_5-A_5|^3+|\alpha_6-A_6|^3\right)^{\dfrac{5}{3}},\\
& \mathrm{fourth~summand:} \left(|\alpha_7-A_7|^{11}+|\alpha_8-A_8|^{11}\right)^{\dfrac{5}{11}},\\
& \mathrm{fifth~summand:} \left(|\alpha_9-A_9|^9+|\alpha_{10}-A_{10}|^9\right)^{\dfrac{5}{9}}.
\end{align*}
We assume $\mathbb{X}_1=\ell_{p_1}^1$ for any $p_1\in (1, \infty)$, $\mathbb{X}_2=\ell_{7}^2$, $\mathbb{X}_3=\ell_{3}^3$, $\mathbb{X}_4=\ell_{11}^2$, $\mathbb{X}_5=\ell_{9}^2.$ Finally, let
$\mathbb{X}=\bigoplus\limits_{i=1}^n{_5~\mathbb{X}_i}.$ \\

\noindent Step II:
Let 
\begin{align*}
& \mathbf{y}_1 =(-9,1,1,3,6,8,0,7,9,12), & \mathbf{y}_2=(7,0,5,9,5,3,2,7,1,\frac{13}{2}),\\
& \mathbf{y}_3=(9,-4,3,-4,7,8,8,1,9,13), & \mathbf{y}_4=(9,1,6,3,-1,-1,-7,0,5,-\frac{3}{2}),\\
& \mathbf{y}_5=(8,3,4,8,6,0,3,5,-3,\frac{5}{2}), & \mathbf{y}_6=(9,3,-8,2,1,2,0,2,7,\frac{11}{2}),\\
& \mathbf{y}_7=(6,-2,9,5,8,1,4,1,5, \frac{11}{2}), & \mathbf{y}_8=(7,6,7,9,8,3,2,1,3, \frac{9}{2}),\\
& \mathbf{y}_9=(-8,0,0,1,-6,9,0,4,9,11), & \mathbf{y}_{10}=(8,9,2,7,5,5,6,9,8,14).
\end{align*}
Define $T:\mathbb{R}^{10}\to \mathbb{R}^{10}$ in such a way that
$$[T]=\left[\mathbf{y}_1^t~ \mathbf{y}_2^t~ \mathbf{y}_3^t~\mathbf{y}_4^t~\mathbf{y}_5^t ~\mathbf{y}_6^t~\mathbf{y}_7^t~\mathbf{y}_8^t~\mathbf{y}_9^t ~\mathbf{y}_{10}^t\right ].$$
It is not difficult to see that $\mathcal{N}(T^*)= span\{(0,0,0,0,0,1,1,1,1,-2)\}$.\\

\noindent Step III: Clearly, the given problem is equivalent to finding the minimum of the collection: 
$$\{\|\mathbf{x_0}-y\|_\mathbb{X}: y\in \mathcal{R}(T)\}.$$
Now, using Theorem \ref{Inequalities}, it is easy to see that
$$\min\{\|\mathbf{x_0}-y\|_\mathbb{X}: y\in \mathcal{R}(T)\}= \left\{1+2^{^{\dfrac{25}{22}}}+\left(1+2^{^{\dfrac{9}{8}}}\right)^{\dfrac{10}{9}}\right\}^{-\dfrac{4}{5}}|\alpha_6+\alpha_7+\alpha_8+\alpha_9-2\alpha_{10}|.$$
Therefore, we have obtained the complete solution to the minimization problem (\ref{Minimization problem}).
\end{example}

\bibliographystyle{amsplain}

\end{document}